\documentclass[reqno,11pt]{amsart} % Hier ggf. 12pt in "[]"  
\usepackage{amssymb}%, varbbb2e} 

\newtheorem{theorem}{Theorem}[section] 
\newtheorem{lemma}[theorem]{Lemma}

\theoremstyle{definition}

\theoremstyle{remark} 
 
\newtheorem{notation}[theorem]{Notation} 
 
\numberwithin{equation}{section} 
 
\def\R{{\mathbb{R}}}
\def\uu{{\underline{u}}}
\def\graph{{\rm graph}\,}
\def\xx#1{{\frac{x_{#1}}{|x|}}}
\def\tang#1#2{{\left(\delta_{#1#2}-\xx #1\xx #2\right)}}
\def\gradeins#1#2{{\left(1+{|D#1|}^2\right)^{#2}}}
\def\tru{{{\rm tr}\,u^{ij}}}
\def\trutt{{{\rm tr}\,u^{\tau\tau'}}}
\def\unn{{u^{\nu\nu}}}

\newcommand{\spc}{\;\;\;\;\;\;\;}
\def\intK{{\stackrel{\circ}{K}}}

\begin{document}
\title{Hypersurfaces of prescribed Gau{\ss} Curvature in
  exterior Domains}

\author{Felix Finster} 
%\address{Max Planck Institute for Mathematics in the Sciences, 
%  Inselstr.\ 22-26, 04103 Leipzig, Germany} 
%\email{Felix.Finster@mis.mpg.de} 
%\thanks{} 

\author{Oliver C.\ Schn\"urer} 
\address{Max Planck Institute for Mathematics in the Sciences,  
  Inselstr.\ 22-26, 04103 Leipzig, Germany} 
%\curraddr{} 
\email{Felix.Finster@mis.mpg.de, Oliver.Schnuerer@mis.mpg.de} 
%\thanks{} 
 
%    General info 
 % since 1980 
 
\date{submitted February 2001, in revised form May 2001.} 
 
%\dedicatory{} 
 
 
\begin{abstract} 
We prove an existence theorem for convex hypersurfaces of prescribed Gau{\ss}
curvature in the complement of a compact set in Euclidean space which are
close to a cone.
\end{abstract} 
 
\maketitle 

\section{Introduction}
\markboth{PRESCRIBED GAUSS CURVATURE IN EXTERIOR DOMAINS}
{F.\ FINSTER, O.\ C.\ SCHN\"URER}
Hypersurfaces of prescribed curvature have been subject to intensive studies. 
To mention a few examples, compact convex hypersurfaces are considered in the
closed case in~\cite{cg}, for Neumann boundary conditions 
in~\cite{ltu}, and for
Dirichlet boundary conditions in~\cite{guan, it, tn, Minkowski 
space}, based on
the method of~\cite{tr}.  In the non-compact case, complete 
convex hypersurfaces
are found in~\cite{cw, pg}, whereas~\cite{kt, iwk} deal with equations of mean
curvature type in exterior domains.  In this paper, we consider the Dirichlet
problem for convex hypersurfaces of prescribed Gau{\ss} curvature in exterior
domains.  More precisely, let 
$\emptyset\not=K \subset \R^n$, $n\ge 2$, be a compact set 
whose boundary $\partial K$ is a smooth submanifold of $\R^n$.  
Then for $0<f \in C^\infty\!\left(\left(\R^n \setminus \intK\right)\times
\R\right)$ and
$u_0 \in C^\infty(\partial K)$ we consider strictly convex hypersurfaces of
prescribed Gauss curvature represented as graphs over 
$\R^n \setminus \intK$, i.~e.\ functions $u\in C^\infty\!\left(
\R^n\setminus\intK\right)$ such that
\begin{equation} \label{exterior Dirichlet problem}
\left\{\begin{array}{rcll}
{\mathcal{K}}[u] \equiv \displaystyle
\frac{\displaystyle\det D^2 u}{\displaystyle
\left(1+ |Du|^2 \right)^{\frac{n+2}{2}}}&\!\!=\!\!&f(x,u)
&\mbox{in~}\R^n\setminus K,\\
u&\!\!=\!\!&u_0&\mbox{on~}\partial K\;.
\end{array}\right. 
\end{equation}
The hypersurfaces should be {\em{close to a cone}} in the sense that
\begin{equation}
\sup|u - r| \;<\; \infty
    \label{eq:cone}
\end{equation}
with $u=u(x)$ and $r \equiv |x|$.  A
function $\underline{u} \in C^\infty\!\left(\R^n \setminus \intK\right)$ 
is called a
{\em{subsolution}} if $\uu$ is strictly convex and
\begin{equation} \label{subsol}
\left\{\begin{array}{rcll}
{\mathcal{K}}[\underline{u}]&\!\!\geq\!\!&f(x,\underline{u})
&\mbox{in~}\R^n\setminus K,\\
\underline{u}&\!\!=\!\!&u_0&\mbox{on~}\partial K \;.
\end{array}\right.
\end{equation}

We prove the following main theorem.
\begin{theorem}\label{main theorem}
Let $u_0 \in C^\infty(\partial K)$ and 
$0<f\in C^\infty\!\left(\left(\R^n \setminus 
\intK\right)\times\R\right)$ with
\begin{equation}
f_z\ge 0,\quad \sup \left(f\cdot r^{n+1}\right)+ 
\sup \left(
\frac{|Df|+|D^2f|}{f} \right)\;<\; 
\infty \;.\label{eq:12a}
\end{equation}
Suppose that $\underline{u}\in C^\infty\!\left(\R^n\setminus
\intK\right)$ is a subsolution~(\ref{subsol}) 
which is close to a cone and satisfies the following decay 
conditions at infinity,
\begin{equation}
|D( \underline{u} - r)| \;=\; {\mathcal{O}}\!
\left( \frac{1}{r} \right) \;,\spc
|D^2( \underline{u} - r)| + |D^3\underline{u}| 
\;=\; {\mathcal{O}}\!\left( \frac{1}{r^2} \right) . \label{eq:regular}
\end{equation}
Then there exists a smooth, strictly convex 
hypersurface of prescribed Gau{\ss} 
curvature~(\ref{exterior Dirichlet problem}) which is close to a cone%
~(\ref{eq:cone}).
\end{theorem}

Note that $u-r$ has no prescribed values at infinity,
because the property that $u$ is close to a cone determines
the asymptotics only up to a bounded function. In particular,
the function $|u-\uu|$ will in general not become small near
infinity. Due to these weak assumptions on the asymptotics,
we cannot expect uniqueness of our solutions. 

We remark that the upper bound for $f$ 
in~(\ref{eq:12a}) is natural because the decay 
conditions for $\uu$, (\ref{eq:regular}), imply that
$f(x,\uu)={\mathcal O}(r^{-n-1})$. Namely, if $\nu$ and $\tau, \tau'$ 
denote unit vectors in radial and axial
directions, respectively, the derivatives of $\underline{u}$ 
decay at infinity
like $|D_{\nu \nu}\underline{u}|, |D_{\nu \tau}\underline{u}| = 
{\mathcal{O}}(r^{-2})$,
$|D_{\tau \tau'}\underline{u}| = {\mathcal{O}}(r^{-1})$. Thus $\det
D^2\underline{u} = {\mathcal{O}}(r^{-n-1})$, and  the differential
inequality~(\ref{subsol}) yields the claimed 
decay rate for $f(x,\uu)$.

The above theorem has the disadvantage that it involves 
the existence of a strictly 
convex subsolution. Therefore, we give methods to find 
subsolutions. In the special case when $K=\overline{B_{\rho_1}}
\equiv \overline{B_{\rho_1}(0)}$ is a ball, the boundary
values are zero, and $f$ decays even faster than in~(\ref{eq:12a}),
we construct an explicit subsolution to get the following result. 
\begin{theorem} \label{thm2}
Let $\rho_1>0$ and $a>2$. If $f \in C^\infty
((\R^n \setminus B_{\rho_1})\times\R)$ 
such that $f_z\ge0$,  
\[ \sup \left(
\frac{|Df|+|D^2f|}{f} \right) \;<\; \infty \;,\spc
0 \;<\; f \;\leq\; (a-1)\: 2^{-\frac{3n}{2}-1}\: 
\rho_1^{a-1}\: r^{1-n-a} \;, \]
then there exists a smooth, strictly convex hypersurface 
of prescribed Gau{\ss} 
curvature,
\[ \left\{\begin{array}{rcll}
{\mathcal{K}}[u] &\!\!=\!\!&f(x,u)
&\mbox{in~}\R^n\setminus B_{\rho_1},\\
u&\!\!=\!\!&0&\mbox{on~}\partial B_{\rho_1}\;,
\end{array}\right. \]
which is close to a cone.
\end{theorem}

A more general method for constructing a subsolution is to 
combine a subsolution
in a bounded set with a subsolution near infinity, 
in the spirit of viscosity solutions.
\begin{theorem} \label{thm3}
Let $0<f\in C^\infty\!\left(\left(\R^n\setminus
\intK\right)\times\R\right)$.
Assume that $\R^n \setminus K=\Omega_1\cup\Omega_2$,
where $\Omega_1$ is a bounded open set, 
$\partial K\subset\partial\Omega_1$, and $\Omega_2$ 
an unbounded open set.
Suppose that $u_0$ is a
strictly convex subsolution in $\Omega_1$,
\[ {\mathcal{K}}[u_0] \;\geq\; f(x,u_0) \spc {\mbox{in $\Omega_1$}}, \]
and that $\underline{u}$ is a strictly convex 
subsolution in $\Omega_2$, both smooth up to the boundary,
without prescribed boundary values.
Furthermore, we assume that
\begin{eqnarray*}
u_0 &<& \underline{u} \spc \mbox{on~} (\R^n \setminus K) \cap 
\partial \Omega_1, \\
\underline{u} &<& u_0 \spc\!\! \mbox{on~}  
\partial \Omega_2\;,
\end{eqnarray*}
that $\uu$ is close to a cone, and (\ref{eq:12a}), 
(\ref{eq:regular}) are satisfied.
Then there exists a smooth, strictly convex hypersurface of prescribed
Gau{\ss} curvature~(\ref{exterior Dirichlet problem}) 
which is close to a cone.
\end{theorem}
This theorem can be applied to the situation when we have a subsolution
with correct boundary values $u_0$ near $\partial K$ and
a subsolution $\uu$ near infinity, which can be glued together
as made precise above.
For $\uu$ one can for example take the subsolution used in the
proof of Theorem~\ref{thm2}, see Lemma~\ref{prthm2}.

The proof of Theorem~\ref{main theorem} is outlined as follows.  It is well
known from \cite{guan} (see also \cite{Minkowski space}) that the Dirichlet
problem for the prescribed Gau{\ss} curvature equation has a strictly convex
solution on smooth
compact domains for smooth data, provided that for the same boundary
values there exists a smooth, strictly convex subsolution $\uu$.
Thus we choose $R_0$ such that $K \subset B_{R_0}$ and consider for any $R>4 
R_0$ the Dirichlet problem
\begin{equation}\label{compact Dirichlet problem}
\left\{\begin{array}{rcll}
{\mathcal{K}}\!\left[u^R\right] &\!\!=\!\!&f\left(x, u^R\right) &
\mbox{in~}B_R\setminus K,\\[.5em]
u^R&\!\!=\!\!&\uu&\mbox{on~}\partial B_R\cup
\partial K.\end{array}\right.
\end{equation}
Our main task is to show that the $C^2$-norms of $u^R$ 
are uniformly bounded in 
$R$. Namely, once this is established, standard Krylov/Shafanov and Schauder
regularity theory yields locally uniform bounds for $u^R$ in any
$C^k$-norm~\cite{T3}.
Using a diagonal sequence argument, we get a subsequence
$u^{R_l}$, $R_l\to\infty$, that converges locally smoothly to a
strictly convex
solution of our original problem. To obtain the required $C^2$-estimates, we 
first use barriers to uniformly bound $u^R - |x|$. These $C^0$-% 
estimates and the convexity of the hypersurfaces allow us to control the 
asymptotic behavior as $R \rightarrow \infty$ of the first derivatives of 
$u^R(x)$ in the region $R/2 \leq |x| \leq R$. Using the 
inner maximum principle
for the largest eigenvalue of $u_{ij}$ as well as standard a priori
estimates for $D^2 u^R$ on $\partial K$, the $C^2$-estimates 
are reduced to controlling $D^2 u^R$ on the outer boundary. We bound the 
tangential second derivatives from above and below by differentiating the 
Dirichlet boundary conditions and using the $C^1$-estimates. By adapting the 
barrier constructions of~\cite{guan, Minkowski space} with constants that have
an appropriate scaling in $R$, we obtain estimates for the 
mixed second derivatives. In contrast to the compact 
case~\cite{guan, Minkowski 
space}, the estimates mentioned so far are so strong 
that when putting them into 
the equation, we immediately get estimates for the 
second derivatives in normal 
direction.

Theorem~\ref{thm2} follows immediately from 
Theorem~\ref{main theorem} by choosing
an explicit subsolution $\underline{u}$, given in Lemma~\ref{prthm2}. To prove
Theorem~\ref{thm3}, we apply Lemma~\ref{prthm3} to get the solution to the
compact Dirichlet problem~(\ref{compact Dirichlet problem}) and then proceed
exactly as in the proof of Theorem~\ref{main theorem}.

\section{Preliminary a priori Estimates}
\begin{notation}
In what follows, lower indices 
denote partial derivatives in $\R^n$. 
We denote the partial derivative of $f$ with respect to the
second argument by $f_z$. We set $(u^{ij})
=(u_{ij})^{-1}$. All other indices are raised with 
respect to the Euclidean metric.
Unless otherwise stated, we will use the Einstein summation 
convention. By $c$
we denote a constant independent of $R$ which may 
change its value from line
to line throughout the text.
\end{notation}

Without loss of generality we can assume 
that $0 \in \stackrel{\circ}{K}$.
We choose a constant 
$L> \max\limits_{\partial K} (u_0 - |x|)$ and set
\begin{equation}
    \overline{u} \;=\; |x|+L \;.
    \label{eq:olu}
\end{equation}
The functions $\overline{u}$ and $\underline{u}$ (as in Theorem~\ref{main 
theorem}) will be used as upper and lower barriers, respectively.
Let $u^R$ be a solution of the Dirichlet 
problem~(\ref{compact Dirichlet problem}).
\begin{lemma}\label{C0 estimate lemma}
As $R \rightarrow \infty$, the functions $u^R$ converge locally uniformly to a
continuous function $u$. 
Moreover, $\uu\le u^R\le \overline{u}$ in $B_R\setminus \intK$.
\end{lemma}
\begin{proof}
 From the maximum principle we deduce that
\begin{equation*}
\uu\le u^R\le \overline{u} \quad\mbox{in~}B_R\setminus K.
\end{equation*}
Hence for $R_1<R_2$,
\begin{equation*}
u^{R_1}\le u^{R_2}\quad\mbox{on~}\partial B_{R_1}\cup\partial K,
\end{equation*}
and again from the maximum principle,
\begin{equation*}
u^{R_1}\le u^{R_2}\quad\mbox{in~}B_{R_1}\setminus K.
\end{equation*}
We conclude that the $u^R$ are monotone in $R$. Their pointwise
limit is convex and thus continuous. So they converge locally
uniformly according to Dini's theorem.
\end{proof}

 From now on we omit the index $R$ and assume that $u$ is
a solution of (\ref{compact Dirichlet problem}) with $R$
fixed sufficiently large.

\begin{lemma}\label{C1 estimate lemma}
$|\nabla u|$ is a priori bounded, uniformly in $R$.
\end{lemma}
\begin{proof}
Since $u$ is strictly convex, $|\nabla u|$ attains its maximum at the
boundary. Tangential derivatives are uniformly bounded there in
view of the Dirichlet boundary conditions. 
The normal derivatives are estimated 
as follows. Let $x \in \partial K$ and $\nu$ the outer unit normal to
$\partial K$ at $x$.
We choose $\lambda_0>0$ independent of $R$ such that the line segment 
$\{x+\lambda \nu,\: 0 \leq \lambda \leq \lambda_0\}$ is contained in $B_{R_0} 
\setminus K$. Using the convexity of $u$ as well as the fact that 
$\underline{u}$ lies below $u$ and $\underline{u}(x)=u(x)$,
\[ \nabla_\nu \underline{u}(x) \;\leq\; \nabla_\nu u(x) \;\leq\; 
\frac{u(x+\lambda_0 \nu) - u(x)}{\lambda_0} \spc ({\mbox{for $x \in \partial
K$)}} \]
($\nabla_\nu u$ denotes the directional derivative).
For $x \in \partial B_R$ and $\nu=x/|x|$, we consider similarly the line 
segment $\{\lambda \nu,\: R_0 \leq \lambda \leq R\}$ and obtain
\[ \frac{u(x) - u(R_0 \nu)}{R-R_0} \;\leq\; \nabla_\nu u(x) \;\leq\; 
\nabla_\nu \underline{u}(x) \spc ({\mbox{for $x \in \partial B_R$)}}. \]
We finally use the uniform $C^0$ bounds of Lemma~\ref{C0 estimate lemma},
in particular that $|u(x)-|x||\le c$.
\end{proof}

The next lemma controls the asymptotic 
behavior of $\nabla u$ as $R$ gets 
large.
\begin{lemma} \label{lemma_C1}
For $\frac{R}{2} \leq |x| \leq R$ let $\nu=x/|x|$ and $\tau$ 
be unit vectors parallel and orthogonal to $x$, respectively.
Then there is a constant $c$ independent of $R$ such that
\begin{eqnarray}
|\nabla_\nu (u-\underline{u})(x)| & \leq & \frac{c}{R} \label{eq:3a} \\
|\nabla_\tau u(x)| & \leq & \frac{c}{\sqrt{R}} \;. \label{eq:3b}
\end{eqnarray}
\end{lemma}
\begin{proof}
Since $u$ is convex and lies above $\underline{u}$ with $u(R \nu) = 
\underline{u}(R \nu)$,
\[ \frac{u(x) - u(R_0 \nu)}{|x| - R_0} \;\leq\; \nabla_\nu u(x) 
\;\leq\; \nabla_\nu \underline{u}(R \nu) \;, \]
and thus
\[ \frac{u(x) - |x|}{|x|-R_0} - \frac{u(R_0 \nu) - R_0}{|x|-R_0} \;\leq\;
\nabla_\nu u(x) - 1 \;\leq\; \nabla_\nu \underline{u}(R \nu) - 1 \;. \]
The $C^0$ estimates of Lemma~\ref{C0 estimate lemma} 
imply that $|u(x)-|x|| \leq
c$, and so the left hand side is ${\mathcal{O}}(R^{-1})$. According
to~(\ref{eq:regular}), $|\nabla_\nu \underline{u}(x) - 1| = 
{\mathcal{O}}(R^{-1})$,
and thus
\[ |\nabla_\nu (u-\underline{u})| \;=\; |\nabla_\nu u(x) - 1| + {\mathcal{O}}
\!\left( \frac{1}{R} \right) \;\leq\; \frac{c}{R}\;. \]
In order to derive~(\ref{eq:3b}), we consider $u$ and the barrier functions 
along the line segment $\{x + \lambda \tau \}$ parametrized by $\lambda \in 
[-\lambda_0, \lambda_0]$, $\lambda_0 = \sqrt{R^2 - |x|^2}$. 
The boundary values of $u$ are $u(\pm \lambda_0) = \uu(\pm \lambda_0)$. Thus 
using that $u$ lies above $\uu$ and is convex, we obtain the estimate
\[ \uu'(-\lambda_0) \;\leq\; u'(-\lambda_0) \;\leq\; u'(\lambda=0) \;\leq\;
u'(\lambda_0) \;\leq\; \uu'(\lambda_0) \;, \]
and thus
\[ | \nabla_\tau u(x) | \;=\;|u'(\lambda=0)|
\;\leq\; \max\left\{|\uu'(\lambda_0)|,
\,|\uu'(-\lambda_0)|\right\} \;. \]
Using~(\ref{eq:regular}),
\[ \underline{u}'(\pm \lambda_0) = 
\left.\nabla_\tau \underline{u}(x)\right|_{x=y} =
\left.\nabla_\tau( \underline{u}(x) - |x|)\right|_{x=y} + 
\left.\nabla_\tau |x|\right|_{x=y} = 
{\mathcal{O}} \!
\left(\frac{1}{R} \right) \pm \frac{\lambda_0}{R}, \]
where we have set $y=x \pm \lambda_0 \tau$.
These estimates imply~(\ref{eq:3b}) provided that $x$ 
is sufficiently close to 
the outer boundary, more precisely if $x$ lies in the strip
\[ R^2 - |x|^2 \;=\; \lambda_0^2 \;\leq\; \kappa R \;, \]
where $\kappa$ is some constant independent of $R$. Now suppose that $x$ lies
outside this strip, $\lambda_0^2 > \kappa R$. We construct straight lines
through $(0, u(x)) \in \R^2$ which are tangential to the hyperbola
$(\lambda, \overline{u}(\lambda)) = (\lambda, \sqrt{|x|^2 + \lambda^2} + L)$.
A short explicit calculation shows that there are exactly two such lines, and 
that they go through the points $(\pm \lambda_1, \overline{u}(\lambda))$ with
\[ \lambda_1^2 \;=\; \frac{|x|^4}{(u(x)-L)^2} - |x|^2\;. \]
Using that $u-|x|$ is bounded uniformly in $R$, one sees that $\lambda_1^2 = 
{\mathcal{O}}(|x|)$, and so we can by increasing $\kappa$ arrange that
\begin{equation}
\lambda_1^2 \;\leq\; \kappa R  \;.
    \label{eq:3c}
\end{equation}
Thus $\lambda_1^2<\lambda_0^2$, meaning that the tangential
points are inside the ball $B_R$.
Since $u$ is convex and lies below $\overline{u}$, the line segments
joining the points $(\lambda=0, u(\lambda=0))$ and $(\pm \lambda_1,
\overline{u}(\pm \lambda_1))$, respectively, both lie above $u$.
Moreover, these line segments are by construction tangential to the
hyperbola $(\lambda, \overline{u}(\lambda))$ at $\lambda=\pm \lambda_1$.
Hence
\[ \overline{u}'(-\lambda_1) \;\leq\;
u'(\lambda=0) \;\leq\; \overline{u}'(\lambda_1) \;. \]
Using~(\ref{eq:olu}), we obtain that
\[ |\nabla_\tau u(x)| \;\leq\; 2\: \frac{\lambda_1}{R}\;, \]
and~(\ref{eq:3c}) gives~(\ref{eq:3b}).
\end{proof}

It remains to derive the $C^2$ a priori estimates. 
As in~\cite{guan}, one obtains that
the second derivatives on $\partial K$ are bounded uniformly in $R$. 
Furthermore, it is well known (see e.\ g.\ \cite{gt}) that considering 
$$w\;=\;\frac{\beta}{2}|Du|^2+\log u_{\xi \xi} $$
in an interior maximum taken over $(x, \xi) \in (\overline{B_R
\setminus K}) \times S^{n-1}$
yields that
\begin{equation*}
\max\limits_{B_R\setminus K}|D^2u|\;\le\; c+
\max\limits_{\partial B_R\cup\partial K}|D^2u|\;.
\end{equation*}
Hence it suffices to bound $|D^2u|$ on the outer boundary $\partial B_R$.
It is here where we shall use 
the assumption $\sup(|Df|+|D^2f|)/f<\infty$.

The following lemma gives an estimate for the tangential 
second derivatives. The 
mixed second derivatives will be treated in the next section, and the normal 
second derivatives in Section~\ref{sec5}.

\begin{lemma}[tangential second derivatives at the outer boundary]
\label{tangential C2 estimates lemma}
\rule{0mm}{1mm} \\
Let $x_0\in\partial B_R$ and $\tau_1$, $\tau_2$ be tangential directions
at $x_0$. Then we have at $x_0$,
\begin{equation*}
\left|u_{\tau_1 \tau_2} - |x|_{\tau_1 \tau_2} \right|
\;\le\; \frac{c}{R^2} \;.
\end{equation*}
\end{lemma}
\begin{proof}
We may assume that $x_0=R\cdot e_n\equiv R\cdot(0,\ldots,0,1)$.
Then
$\partial B_R$ is represented locally as $\graph\omega$, where
\begin{equation*}
\begin{array}{rcl}
\omega:\hat B_R\equiv\{\hat x\in\R^{n-1}:|\hat x|<R\}&\to&\R,\\
\hat x&\mapsto&\sqrt{R^2-|\hat x|^2}.
\end{array}
\end{equation*}
According to the Dirichlet boundary conditions,
\begin{equation*}
(u-\uu)(\hat x,\omega(\hat x))=0.
\end{equation*}
We differentiate
twice with respect to $\hat x^i$, $\hat x^j$, 
$1\le i,\,j\le n-1$ and obtain that
at $x_0$,
\begin{equation*}
(u-\uu)_{ij}+(u-\uu)_n\: \omega_{ij}=0 \:.
\end{equation*}
According to the decay conditions at infinity~%
(\ref{eq:regular}),
\begin{equation*}
\left|\uu_{ij} - |x|_{ij} \right| \;=\; 
{\mathcal{O}}\!\left( \frac{1}{R^2} \right) , 
\end{equation*}
and furthermore
\begin{equation*}
\omega_{ij}(x_0) \;=\; -\frac{\delta_{ij}}{R} \;.
\end{equation*}
Thus the result follows in view of Lemma~\ref{lemma_C1}.
\end{proof}

\section{Mixed $C^2$-Estimates} \label{sec4}

\begin{lemma}[Mixed second derivatives at the outer boundary]
\label{Mixed second derivatives}
 For $x_0 \in \partial B_R$ let $\tau$ and $\nu$ 
be unit vectors in tangential 
and normal directions, respectively. Then
\begin{equation}
|u_{\tau\nu}|(x_0)\le\frac{c}{\sqrt{R}} \:. \label{eq:msd}
\end{equation}
\end{lemma}
The proof of this Lemma is split up into several lemmata.
\par
As in the proof of Lemma \ref{tangential C2 estimates lemma} we
may assume that $x_0=R\cdot e_n$ and represent $\partial B_R$
locally as $\graph\omega$ with $\omega(\hat x)=\sqrt{R^2-|\hat x|^2}$.
We take the logarithm of the differential equation in 
(\ref{compact Dirichlet problem}),
$$\log\det u_{ij}-\frac{n+2}{2}\log
\left(1+|Du|^2\right)=\log f(x,u),$$
and differentiate with respect to $x^k$,
\begin{equation}\label{differentiated equation}
u^{ij}u_{ijk}-\frac{n+2}{1+|Du|^2}
u^iu_{ik}=\frac{f_k+f_zu_k}{f}.
\end{equation}
This gives
a motivation for introducing the linear differential operator $L$ by
$$Lw:=u^{ij}w_{ij}-\frac{n+2}{1+|Du|^2}u^iw_i.$$
Furthermore, we define for $t<n$ the linear operator
$$T:=\frac{\partial}{\partial x^t}+\omega_{tr}(0)\:x^r
\frac{\partial}{\partial x^n}\equiv\frac{\partial}{\partial x^t}
-\frac{x_t}{R}\frac{\partial}{\partial x^n},$$
where we used the convention $\omega_{tn}=0$ (thus we sum over
$r=1,\ldots,n-1$). In what follows we restrict attention to the domain
$\Omega_\delta:=B_\delta(x_0)\cap B_R$ with $\delta\le\frac{R}{2}$,
notice that $\Omega_\delta\subset B_R\setminus K$.

\begin{lemma}\label{LTu-u estimates}
The function $u-\underline{u}$ satisfies the following estimates,
$$\begin{array}{rcll}
|T(u-\underline{u})|&\le& \displaystyle \frac{c}{\sqrt{R}}&\mbox{in~}
\Omega_\delta, \\[.9em]
|T(u-\underline{u})|&\le& \displaystyle \frac{c}{R^2}\cdot|x-x_0|^2&\mbox{on~}
\partial B_R, \\[.9em]
|LT(u-\underline u)|&\le& \displaystyle 
c+\frac{c}{R^2}\:\tru&
\mbox{in~}\Omega_\delta \;,
\end{array}$$
where $\tru\equiv u^{ij}\delta_{ij}$.
\end{lemma}
\begin{proof}
Note that $|\omega_i|\le c$, $|\omega_{ij}|\le\frac{c}{R}$
and $|\omega_{ijk}|\le\frac{c}{R^2}$. The first inequality follows directly 
 from the $C^1$ estimates of Lemma~\ref{lemma_C1}, whereas for the second 
inequality we use furthermore that $(u-\uu)_t+(u-\uu)_n\omega_t=0$ on 
$\partial B_R$, the decay properties of $\omega_{ijk}$, and the 
fact that $u-\underline{u}$ vanishes on $\partial B_R$. To prove the last 
inequality, we apply (\ref{differentiated equation}),
the relation $u^{ij}u_{jk}=\delta^i_k$, 
and the $C^1$-estimates to 
obtain that
$$|LT(u-\underline{u})|\le c\cdot\frac{|Df|}{f}+\frac{c}{R}
+c\cdot|D^2\underline{u}|+c\cdot\tru\cdot\left(|D^3\underline{u}|
+\frac{1}{R}|D^2\underline{u}|\right).$$
Now we use the conditions~(\ref{eq:12a}) and (\ref{eq:regular}).
\end{proof} 

The function $\vartheta$ introduced in the next Lemma will be
the main part of a barrier function which we shall
construct in what follows. 
\begin{lemma}\label{small theta}
There exists a positive constant $\varepsilon$ independent
of $R$ such that
$$\vartheta:=(u-\underline{u})+\frac{1}{\sqrt{R}}\: d
-\frac{1}{2R^{\frac{5}{4}}}\: d^2$$
fulfills the estimates
$$\left\{\begin{array}{rcll}L\vartheta&\le&-\varepsilon \:R^{\frac{3}{4n}}
-\varepsilon \:R^{-\frac{5}{4}}\:\tru
&\mbox{in~}\Omega_\delta,\\[.3em]
\vartheta&\ge&0&\mbox{on~}\partial\Omega_\delta \;, \end{array}\right.$$
provided that $\delta = R^{\frac{3}{4}}$ and $R$ is sufficiently 
large. Here $d=R-|x|$ is the distance from $\partial B_R$.
\end{lemma}
\begin{proof}
It is obvious that $\vartheta\ge0$ on $\partial\Omega_\delta$.
We fix $x_0 \in \Omega_\delta$ and set $\nu=x_0/|x_0|$. Let
$\tau, \tau'$ belong to an orthonormal basis for the orthogonal 
complement of
$\nu$ which we choose such that the submatrix $u^{\tau \tau'}$ is 
diagonal. Assume that $\nu$ and $\tau$, $\tau'$ correspond to the
indices $n$ and $1,\,\ldots,\,n-1$, respectively.
We use the Einstein summation convention for $\tau, \tau'$.
The matrix $u^{ij}$ is positive, and thus testing with the vectors
$\nu \pm \tau$ gives
\begin{equation}
|u^{\nu \tau}| \;\leq\; \frac{1}{2}\: (u^{\nu \nu} + u^{\tau \tau}) \;.
\label{mix}
\end{equation}
We introduce the abbreviation
\[ \trutt \;=\; u^{ij}\tang ij . \]
Direct computations using~(\ref{eq:regular}) and~(\ref{mix}) give
\begin{eqnarray*}
L u &=& u^{ij} u_{ij} - (n+2)\: \frac{|Du|^2}{1+|Du|^2} \;\leq\; c, \\
L \underline{u} &\geq& -c \:+\: u^{\tau 
\tau'}\:\underline{u}_{\tau \tau'} \:+\: 2\: 
u^{\tau \nu}\:\underline{u}_{\tau \nu} 
\:+\: u^{\nu \nu}\:\underline{u}_{\nu \nu} \\
&\geq& -c + \left(\frac{1}{|x|} - \frac{c}{|x|^2} \right) \:
\trutt - \frac{c}{|x|^2} \:u^{\nu \nu}, \\
L\vartheta&\le&c\cdot\left(1+\frac{1}{\sqrt{R}}
\right)-\left(\frac{1}{|x|}-\frac{c}{|x|^2}+\frac{1}{|x|}
\left(\frac{1}{\sqrt{R}}-\frac{d}{R^{\frac{5}{4}}}\right)\right)
\trutt\\
&&{}-\left(\frac{1}{R^{\frac{5}{4}}}-\frac{c}{|x|^2}\right)\unn,
\end{eqnarray*}
and thus for $R$ sufficiently large,
$$L \vartheta \;\le\; c-\frac{1}{2R}\trutt-\frac{1}
{2R^{\frac{5}{4}}}\: \unn.$$

Expanding the determinant and using that $u^{\tau \tau'}$ is diagonal gives
\begin{eqnarray*}
\det u^{ij} &=&
\det\left(\begin{array}{ccccc}
u^{11} & 0       & \cdots & 0             & u^{1n}   \\
0      & \ddots  & \ddots & \vdots        & \vdots   \\
\vdots & \ddots  & \ddots & 0             & \vdots   \\
0      & \cdots  & 0      & u^{n-1\, n-1} & u^{n-1\, n}\\
u^{1n} & \cdots  & \cdots & u^{n-1\, n}   & u^{nn}   \\
\end{array}\right) \\[.3em]
&=& \prod_{i} u^{ii} \:-\: \sum_\tau |u^{n \tau}|^2 
\:\prod_{\tau' \neq \tau} u^{\tau' \tau'} \;\leq\; \prod_{i} u^{ii} \;.
\end{eqnarray*}
Hence the inequality for geometric and arithmetic means as well
as~(\ref{compact Dirichlet problem}) and 
(\ref{eq:12a}) show that for large values of $R$,
\begin{eqnarray*}
L\vartheta&\le&c-\frac{1}{c}\cdot(\det u_{ij})^{-\frac{1}{n}}
\cdot R^{-\frac{n-1}{n}-\frac{5}{4}\frac{1}{n}}
-\frac{1}{4 R^{\frac{5}{4}}}\:\tru\\
&\le&c-\frac{1}{c}\: |x|^{\frac{n+1}{n}} 
\:R^{-\frac{1}{n}\left((n-1)+\frac{5}{4}
\right)} -\frac{1}{4 R^{\frac{5}{4}}}\:\tru \\
&\le&-\frac{1}{c}\cdot R^{\frac{3}{4n}} -\frac{1}{4 R^{\frac{5}{4}}}
\:\tru.
\end{eqnarray*}
\end{proof} 

\begin{lemma}\label{big Theta}
There exists a positive constant $A$ independent of $R$ such that
$$\Theta \;:=\; \vartheta+A\cdot\frac{1}{R^2}\cdot|x-x_0|^2\pm T(u-
\underline{u})$$
satisfies the inequalities
\begin{equation}
    \left\{\begin{array}{rcll}L\Theta&\le&0&\mbox{in~}\Omega_\delta,\\
    \Theta&\ge&0&\mbox{on~}\partial\Omega_\delta,\end{array}\right.
    \label{eq:34}
\end{equation}
where $\delta=R^{\frac{3}{4}}$ and $\vartheta$ is as in Lemma
\ref{small theta}.
\end{lemma}
\begin{proof}
According to Lemma~\ref{small theta}, the condition
$\Theta\ge0$ on $\partial\Omega_\delta$ follows if 
$$A\cdot\frac{1}{R^2}\cdot|x-x_0|^2\pm T(u-\underline{u})\ge0
\quad\mbox{on~}\partial\Omega_\delta.$$
In view of Lemma \ref{LTu-u estimates}, this can be arranged by choosing $A$ 
sufficiently large. The property $L\Theta\le 0$ now follows from the
inequality
$$-\varepsilon R^{\frac{3}{4n}}-
\varepsilon R^{-\frac{5}{4}}\: \tru+c\cdot\frac{A}{R}
+c+c\cdot\frac{1+A}{R^2}\:
\tru\le 0,$$
which holds for $R$ sufficiently large.
\end{proof} 

\begin{proof}[Proof of Lemma \ref{Mixed second derivatives}.]
The maximum principle applied to~(\ref{eq:34}) yields that $\Theta \geq 0$ 
in $\Omega_\delta$. Since $\Theta(x_0)=0$, it follows that
$$\Theta_\nu(x_0)\ge0 $$
with $\nu=-\frac{x_0}{|x_0|}$. Thus we obtain
$$\vartheta_\nu(x_0)\ge|(T(u-\underline{u}))_\nu|(x_0) \;, $$
and this finally gives~(\ref{eq:msd}).
\end{proof}

\section{Remaining A Priori Estimates} \label{sec5}
The decay rate in Lemma \ref{Mixed second derivatives} has been
chosen such that the following lemma follows immediately. 

\begin{lemma}[Double normal $C^2$-estimates at the outer boundary]
\rule{0mm}{1mm} \\
Under the assumptions of Lemma~\ref{Mixed second derivatives}
and of one of our Theorems~\ref{main theorem}-\ref{thm3},
\begin{equation*}
|u_{\nu\nu}|(x_0)\le c.
\end{equation*}
\end{lemma}
\begin{proof}
We choose for fixed $x_0 \in \partial B_R$ an orthonormal basis as in
Lemma~\ref{small theta}, now such that the submatrix $u_{\tau \tau'}$
is diagonal. We expand the determinant,
\begin{eqnarray*}
f(x,u)\cdot\gradeins u{\frac{n+2}{2}} &=& \det u_{ij} 
\;=\; u_{nn}\cdot\prod\limits_{i<n}u_{ii}
-\sum\limits_{k<n}u^2_{kn}\cdot\prod\limits_{k\not=i<n}u_{ii} \\
&=& u_{nn}\cdot\prod\limits_{i<n}u_{ii}
-\prod\limits_{i<n}u_{ii}\cdot
\sum\limits_{k<n}u^2_{kn}\:\frac{1}{u_{kk}} \;.
\end{eqnarray*}
Now we substitute in the estimates of 
Lemma~\ref{tangential C2 estimates lemma}
and Lemma~\ref{Mixed second derivatives},
\begin{eqnarray*}
u_{nn}&\le& \frac{c}{\displaystyle \prod\limits_{i<n} u_{ii}}\cdot  
f(x,u)+\sum\limits_{k<n}\frac{u^2_{kn}}{u_{kk}} \\
&\le&c\cdot R^{n-1}\cdot R^{-n-1}
+\sum\limits_{k<n} \frac{\left(\frac{c}{\sqrt{R}}\right)^2}
{\frac{c}{R}} \;\leq\; c \;.
\end{eqnarray*}
\end{proof}

\section{Barrier Constructions}
The following lemma gives a simple example of a barrier construction.
\begin{lemma} \label{prthm2}
Let $K=\overline{B_{\rho_1}}(0)$ with $\rho_1>0$ and $u_0 \equiv 0$.
Suppose that $f \in C^\infty \!\left(\R^n
\setminus \stackrel{\circ}{K} \right)$
satisfies
\[ 0 \;<\; f \;\leq\; (a-1)\: 2^{-\frac{3n}{2}-1}\: 
\rho_1^{a-1}\: r^{1-n-a} \]
with $a>2$. Then there is a strictly convex subsolution $\underline{u}$ which
is close to a cone such that~(\ref{eq:regular}) is satisfied.
\end{lemma}
\begin{proof}
We introduce functions $\varphi$ and $\psi$ by
\begin{eqnarray*}
\varphi(r) &=& \frac{a-1}{2}\: \rho_1^{a-1}\: r^{-a} \\
\psi(r) &=& -\int_{\rho_1}^r \left( \int_\rho^\infty \varphi(\tau)\: d\tau 
\right) d\rho
\end{eqnarray*}
and define $\underline{u}$ by
\[ \underline{u} \;:\; \R^n \setminus B_{\rho_1} \rightarrow \R \;:\;
x \mapsto |x| - \rho_1 + \psi(|x|) \;. \]
This $\underline{u}$ is zero on $\partial B_{\rho_1}$, 
is close to a cone, and 
satisfies the regularity conditions~(\ref{eq:regular}). Furthermore,
\[ 0 \;\leq\; - \psi'(r) \;=\; \int_r^\infty \varphi(\tau)\: d\tau \;\leq\;
\int_{\rho_1}^\infty \varphi(\tau)\: d\tau \;=\; \frac{1}{2}\;. \]
We compute the Gau{\ss} curvature of $\graph\uu$,
\begin{eqnarray*}
\uu_i&=& \xx i+\psi'\xx i,\quad
|D\uu|=1+\psi', \\
\uu_{ij}&=& \frac{1}{|x|}(1+\psi')\tang ij
+\varphi(|x|)\xx i\xx j \;.
\end{eqnarray*}
We conclude that in $\R^n\setminus B_{\rho_1}$
\begin{eqnarray*}
{\mathcal{K}}[\uu]&=&\varphi\cdot r^{1-n}\cdot(1+\psi')^{n-1}\cdot
\left(1+(1+\psi')^2\right)^{-\frac{n+2}{2}}\\
&\ge&(a-1)\cdot\rho_1^{a-1}\cdot r^{1-n-a}\cdot 2^{-n}
\cdot 2^{-\frac{n+2}{2}}.
\end{eqnarray*}
\end{proof}

The above construction
 could be extended to more general domains and non-zero boundary
values, provided that the boundary values are sufficiently close to a cone,
i.\ e.\ $|u_0-|x||_{C^2} < \varepsilon$, and $f$ is sufficiently small.
For non-zero boundary values, it seems preferable to use the following
construction of a subsolution, based on the idea of viscosity subsolutions.

\begin{lemma} \label{prthm3}
For a smooth domain $\Omega\subset\R^n$ and a smooth positive function
$f:\overline{\Omega}\times\R\to\R$ with $f_z\ge0$, we assume that 
there exist two smooth, strictly convex subsolutions, i.\ e.
$$u_i:\overline{\Omega_i}\to\R, 
\quad {\mathcal{K}}[u_i]\ge f(x,u_i)\mbox{~~in~}
\Omega_i,\quad i=1,2,$$
such that $\Omega_1\cup\Omega_2=\Omega$ $(\Omega_i=\mbox{open})$, 
$u_i<u_{3-i}\mbox{~on~}
\partial\Omega_i$, $i=1,2$, and $\{u_1=u_2\}\Subset\Omega_1\cap
\Omega_2$ (in particular $\{u_1=u_2\}\cap\partial\Omega=
\emptyset$). Then there exists a smooth function $u:\overline{\Omega}
\to\R$ such that
$$\left\{\begin{array}{rcll} {\mathcal{K}}[u]&=&f(x,u)&
\mbox{in~}\Omega,\\[.2em]
u&=&\max\{u_1,u_2\}&\mbox{on~}\partial\Omega,\\[.2em]
u&\ge&\max\{u_1,u_2\}&\mbox{in~}\Omega.\end{array}\right.$$
\end{lemma}
\begin{proof}
We choose a smooth, strictly convex function $w:\overline{\Omega}\to
\R$ such that
$$\left\{\begin{array}{rcll}w&\ge&\max\{u_1,u_2\}&\mbox{in~}
\Omega,\\[.2em] w&=&\max\{u_1,u_2\}&\mbox{on~}\partial\Omega
\end{array}\right.$$
and define $f_w(x):={\mathcal{K}}[w](x)$. 
Then $w$ is a solution of the Dirichlet
problem
$$\left\{\begin{array}{rcll} {\mathcal{K}}[w]&=&f_w(x)&
\mbox{in~}\Omega,\\[.2em]
w&=&\max\{u_1,u_2\}&\mbox{on~}\partial\Omega.\end{array}\right.$$
We may assume that 
$f_w<f(\cdot,\max\{u_1,u_2\})$ in $\Omega$ because
otherwise we solve a Dirichlet problem 
(with small $f$) by using 
\cite{guan} with $w$ as the subsolution, and take the
solution to this Dirichlet problem instead of $w$. 
Finally, we apply the
continuity method and solve the Dirichlet problems
$$\left\{\begin{array}{rcll}
{\mathcal{K}}\left[u^t\right]&=&(1-t)f_w(x)+tf(x,u^t)&
\mbox{in~}\Omega,\\[.5em]
u^t&=&\max\{u_1,u_2\}&\mbox{on~}\partial\Omega,\\[.2em]
u^t&\ge&\max\{u_1,u_2\}&\mbox{in~}\Omega,\end{array}\right.$$
for smooth strictly convex functions $u^t:\overline{\Omega}\to\R$,
$0\le t\le 1$. The solvability of this Dirichlet problem 
is a consequence of 
the a priori estimates of \cite{guan} and the fact that
$$u^t>\max\{u_1,u_2\}\quad\mbox{in~}\Omega\mbox{~~for~}
0\le t<1.$$ 
The last inequality follows from the
maximum principle because $(1-t)f_w(x)+tf(x,\max\{u_1,u_2\})
<f(x,\max\{u_1,u_2\})$ 
for $0\le t<1$. For $t=1$ we obtain the function $u$ we are
looking for.
\end{proof}

\begin{appendix}
\section{Notes}
We remark that the $C^2$-estimates at the outer boundary
$\partial B_R$ cannot be obtained from the estimates on compact
domains by a simple scaling argument. The reason is that if one
considers the function $\frac{1}{R}\: u(Rx)$ near $B_1 \!\setminus\!
\left(\frac{1}{R}\cdot K \right)$, the corresponding rescaled
function $f$ tends to zero near $\partial B_1$ as $R \rightarrow
\infty$.

The purpose of this paper was to find a method for 
proving the existence of strictly convex hypersurfaces of 
prescribed Gau{\ss} curvature in exterior domains. But for 
technical simplicity, we did not consider the most
general situation to which our methods apply. 

First of all,
our regularity assumptions could be weakened; indeed, it is 
sufficient to assume the same regularity 
as in the compact case, i.\ e.\ similar to~\cite{Minkowski space}
\[ \partial K \in C^4 \;,\spc \underline{u} \in 
C^4\!\left(\R^n \setminus \intK\right) \;,\spc
f \in C^2\!\left(\left(\R^n \setminus \intK\right)\times\R\right)\;. \]
Furthermore, it is possible to treat
equations of the form
$$\det u_{ij}=f(x,u,Du)$$
with $f>0$. We do not consider
 this case here because the explicit construction
of a subsolution becomes more complicated for this more general
ansatz.

Moreover, it is possible to remove the technical condition 
$f_z\ge 0$ by using
a mapping degree argument and assuming that $f_z\ge 0$ only outside
the set $\{(x,z)\in(\R^n\setminus K)\times\R:\uu(x)\le z\le
\overline{u}(x)\}$. But the uniqueness as described in Lemma~%
\ref{C0 estimate lemma} is lost, i.\ e.\ the limit of $u^{R_k}$
as $R_k\to\infty$ may depend on the choice of the subsequence.
We also remark that the decay condition
$$\sup\left(
\frac{|Df|+|D^2f|}{f}\right)<\infty$$ 
is clearly only needed between the two barriers.

Finally, one could consider hypersurfaces of prescribed Gau{\ss}
curvature which, instead of being close to a cone, have a 
different asymptotic behavior near infinity. One example 
would be to take a smooth convex domain $\Omega\subset\R^{n-1}$, 
$0\in\Omega$, and to consider hypersurfaces which are close
to the cone $\R_{\,+}\cdot(\partial\Omega\times\{1\})$.
It even seems possible to adapt our methods to hypersurfaces of
prescribed Gau{\ss} curvature which are close to a more
general convex hypersurface, provided that 
its Gau{\ss} curvature decays 
sufficiently fast at infinity, and that 
one can find barriers that give good $C^1$-estimates.
\end{appendix}

\bibliographystyle{amsplain}

\end{document}